\documentstyle[12pt]{article}
\renewcommand{\baselinestretch}{1.3}

\newtheorem {th}{Theorem}[section]
\newtheorem {lem}[th]{Lemma}
\newtheorem {pr}[th]{Proposition}
\newtheorem {cor}[th]{Corollary}
\newtheorem{defn}[th]{Definition}
\newtheorem{defns}[th]{Definitions}
\newtheorem{conj}{conjecture}
\def\Cox{\hfill \Box}

\def\parent{{\theta_1}}
\def\child{{\theta_2}}
\def\dd{\delta}

\def\iid{\mbox{independent, identically distributed}}

\def\ee{\epsilon}

\def\c1{c_1}
\def\d1{d_1}

\def\v0{{\vec 0}}

\def\E{{\bf{E}}}
\def\P{{\bf{P}}}

\def\R{{\bf{R}}}
\def\Z{{\bf{Z}}}
\def\A{{\cal{A}}}

\def\F{{\cal{F}}}

\def\B{{\cal{B}}}

\def\c1{c_1}
\def\d1{\delta_1}

\def\E{{\bf{E}}}

\def\|{\, | \, }

\def\Tree{{\bf T}}

\begin{document}

\begin{titlepage}
\begin{center}
{\large \bf THE CONTACT PROCESS ON TREES} \\
\end{center}

\vspace{2ex}
\begin{quote} \raggedleft
Robin Pemantle \footnote{This research was supported by an NSF 
  postdoctoral fellowship and a Mathematical Sciences Institute
  research fellowship.  Presently at the University of Wisconsin-Madison}\\
Department of Mathematics \\
Cornell University \\
Ithaca, NY 14853
\end{quote}

\vfill

{\bf ABSTRACT:} \break

\renewcommand{\baselinestretch}{1.0}\large\normalsize
\noindent{The} contact process on an infinite homogeneous tree is
shown to exhibit at least two phase transitions as the infection
parameter $\lambda$ is varied.  For small values of $\lambda$ 
a single infection eventually dies out.  For larger $\lambda$ the
infection lives forever with positive probability but eventually
leaves any finite set.  (The survival probability is a continuous 
function of $\lambda$, and the proof of this is much easier than
it is for the contact process on $d$-dimensional integer lattices.)  
For still larger $\lambda$ the infection
converges in distribution to a nontrivial invariant measure.
For an $n$-ary tree, with $n$ large, the first of these transitions 
occurs when $\lambda \approx 1/n$ and the second occurs when 
$1/2\sqrt{n} < \lambda < e/\sqrt{n}$.  Nonhomogeneous trees whose
vertices have degrees varying between $1$ and $n$ behave
essentially as homogeneous $n$-ary trees, provided that vertices
of degree $n$ are not too rare.  In particular, letting $n$ go
to $\infty$, Galton-Watson
trees whose vertices have degree $n$ with probability that does
not decrease exponentially with $n$ may have both phase transitions
occur together at $\lambda = 0$.  The nature of the second phase 
transition is not yet clear and several problems are mentioned in
this regard.

\hspace{1.2in} Galton-Watson tree, periodic tree

\renewcommand{\baselinestretch}{1.5}\large\normalsize
\end{titlepage}

\section{Introduction} \label{intro}

This paper studies the contact process on trees.  A tree is
just a connected unoriented graph without cycles.  A brief
description of the contact process is as follows.  The state
space, in this case a tree, is any undirected graph.  At
any instant in time, some vertices of the graph are infected
and some are healthy.  The infected vertices recover (turn into
healthy vertices) at poisson rate $1$, independently.  Each 
infected vertex independently infects each of its healthy neighbors
at poisson rate $\lambda$ for some fixed parameter $\lambda$.
More complete
descriptions of this are available in [Li] and [Du], along
with a construction (which I will later describe) of the 
process from a poisson process called the graphical representation. 

There are practical and mathematical reasons to study the
contact process on trees.  The practical reason is that the
process was developed to model the spread of infectious disease.  The 
vertices of the graph represent individuals susceptible
to the disease and the edges represent pairs of individuals
who may have frequent contact.  Traditionally, the graph is taken
to be the $n$-dimensional integer lattice, probably because that
makes sense for other interacting particle systems such as the
Ising model for ferromagnetism.  While the ``real'' graph is a
large finite graph, a tree is at least as likely as a lattice to
serve as a local approximation to the real graph for a population.

One mathematical reason for studying the contact process on trees
is to better understand the correspondence between behaviour of
processes on trees and limits of behaviours on lattices as the
dimension of the lattice goes to infinity.  For other processes
such as random walks and the Ising model, there is a critical
dimension beyond which the lattice versions of the process behave
in essential ways like the tree version.   %%% GET REFERENCES !!!

Another reason to study trees is the following observation about
the contact process on different graphs.  If $G$ and $H$ are graphs
with a map $\phi$ from the vertices of $G$ to the vertices of $H$,
such that $(x,y)$ is an edge of $G$ whenever $(\phi x , \phi y)$
is an edge of $H$, then the contact process may be coupled so
that the number of infections on $G$ is always at least the number
of infections on $H$.  (Letting infections and recoveries on $G$
proceed as on $H$ but subject to the restriction that there can
be at most one infection in any $\phi^{-1}(x)$ gives a process
with equal recovery rate but smaller infection rate than the 
contact process on $G$.  This process has as many infections as the
process on $H$ and by monotonicity, has fewer infections than the
contact process on $G$.)  Since any graph has a tree for a universal
cover, this relates the contact process on any graph to the contact
process on a tree.  In particular, among all graphs where every
vertex has the same degree $n+1$, the tree has the best survival 
probability for the contact process, hence the critical value
$\lambda_1$ defined below is least for the tree.

The main reason that the contact process on a tree is worth
studying is that it
exhibits a multiple phase transition.  In this respect the behaviour
differs remarkably from any of the usual interacting particle
systems on the lattice.  In fact a major recent result [BG]
says that no such monkey business can happen with the contact process.
That such behaviour may be the rule for trees, rather than the 
exception, is made plausible by the recent work of Grimmett, Newman
and Wu [GN], [NW], which shows a similar phenomenon for percolation on
a graph which is the integers crossed with a tree.  Since the
contact process on a graph $G$ is very similar to oriented percolation
on $G \times \Z$, it is not surprising that unoriented percolation
on the integers crossed with a tree should behave like the contact
process on a tree.  

The results in this paper are concerned with
the existence and location of the multiple phase transitions. 
The theory of the contact process on a tree has by no means been
tidied up.  In particular, I believe but am far from being able to 
prove that there are only two phase transitions.  The rest of this
introductory section will define the phase transitions and
delimit what I know about their existence.  

Let $\Tree$ be an infinite tree with a distinguished vertex $\rho$
called the root.  Let $2^\Tree$ be the set of subsets of the
vertices of $\Tree$, and let $\A \subseteq 2^\Tree$ be the 
finite subsets of vertices.  The contact process can be viewed
as a random variable taking values in the space of functions
$\xi$ with domain $\R^+$ and range $2^\Tree$, where $\xi (t)$
is the set of vertices infected at time $t$.  If the initial
configuration $\xi (0)$ is almost surely in $\A$, then $\xi (t)$
stays in $\A$ almost surely for all time.  I will not prove
statements on this level, or even formally construct the contact 
process; the reader may easily do this from [Du].  The idea is to 
construct the graphical representation, which is merely a set of
independent poisson processes.  There is a rate 1 poisson
process for each vertex giving the recovery times at that vertex
(a vertex ignores recovery times while it is healthy) and there is
a rate $\lambda$ poisson process for each oriented edge giving
the {\em attempted infection times}, where an attempted infection
from an infected node to a healthy node transmits the infection
but all other attempted infections are ignored.  It is also
possible to construct the process directly from its semigroup [Li]. 
Here is some notation.

For any configuration $\eta \in 2^\Tree$,
use $\P_\eta (\xi (t) \mbox{ has property X})$ to denote 
probabilities for the contact process started from the initial
condition $\xi (0) = \eta$.  Let $\P_v$ denote $\P_{ \{ v \} }$
in this regard.  A measure $\mu$ is invariant for the contact
process if $\xi (t)$ has law $\mu$ for any $t$ whenever
$\xi (0)$ has law $\mu$.  Since the contact process is an attractive
system, it is easily seen [Li] that there is a stochastically
greatest invariant measure ${\overline \mu}$ that can be gotten
by taking the limit of the laws of $\xi (t)$ as $t \rightarrow \infty$
when $\xi (0)$ is the configuration with every vertex infected.
(The contact process, since it is an attractive system, is monotone
in the following sense. 
Two contact processes run from initial states where the
infections of one are a subset of the infections of
the other may be coupled so as always to remain in this relation.
This and its consequences will be used freely in what follows, and
sometimes without explicit mention.)
The following are all possible definitions of critical values for
the infection parameter $\lambda$.

\begin{defns}
\begin{eqnarray*}
\lambda_1 & = & \inf \{\lambda \,:\, \P_\rho (\xi (t) \neq \emptyset
    \mbox{ for all } t) > 0 \} \\[2ex]
\lambda_a & = & \inf \{\lambda \,:\, \P_\rho (\rho \in \xi (t)
    \mbox{ for arbitrarily large }t ) > 0 \} \\[2ex] 
\lambda_b & = & \inf \{\lambda \,:\, \lim\sup_{t \rightarrow \infty} 
    \P_\rho (\rho \in \xi (t)) > 0 \} \\[2ex]
\lambda_2 & = & \inf \{\lambda \,:\, \lim\inf_{t \rightarrow \infty} 
    \P_\rho (\rho \in \xi (t)) > 0 \} \\[2ex]
\lambda_{c} & = & \inf \{\lambda \,:\, \xi (t) \mbox{ conditioned
    on not being empty } {\stackrel {{\cal D}}
    {\Rightarrow}} {\overline \mu} \neq \dd_0 \}
\end{eqnarray*}
where ${\overline \mu}$ is the upper inveriant measure and $\dd_0$
is the point mass at the empty set.
\end{defns}

In words: $\lambda_1$ is the value at which the process can 
survive forever with nonzero probability, above $\lambda_a$ the
root gets reinfected infinitely often, above $\lambda_b$ the
root is reinfected with probability bounded away from zero at
arbitrarily large times, above $\lambda_2$ these arbitrarily 
large times are all times, and above $\lambda_c$ there is complete
convergence to a nontrivial measure
(which has been shown for one-dimensional,
and more recently any dimensional integer lattice).
It is immediate that $\lambda_1 \leq \lambda_a \leq \lambda_b
\leq \lambda_2 \leq \lambda_c < \infty$, with the last inequality
following from the classical results one the contact process on the
integers, which are a subgraph of any infinite tree.  I believe 
the following relations hold among the critical values.
\begin{conj} \label{cj1}
For any infinite tree, $\lambda_a = \lambda_b = \lambda_2 
= \lambda_c$.
\end{conj}
Call a vertex $v$ of an infinite tree $\Tree$ {\em esential} if
$\Tree \setminus \{ v \}$ has at least two infinite components.
Say that a tree has strongly exponential growth if there is some
$r \geq 1$ such that for any $v \in \Tree$ there are at least
$3$ essential vertices $w$ whose graph-distance from $v$ is $r$.
\begin{conj} \label{cj2}
For any tree with strongly exponential growth,
either all the above critical values are
zero, or else $\lambda_1 \neq \lambda_a$.  
\end{conj}
The nature of the phase transitions other than $\lambda_1$ is
evidently unclear, but regarding the first phase transition I
believe:
\begin{conj} \label{cj3}
For any tree, the survival probability starting from a single
infection is a continuous function of $\lambda$.
\end{conj}
Conjecture~\ref{cj3} is undoubtedly difficult since it includes
the case of the one-dimensional integer lattice, a recent and
difficult result. Theorem~\ref{leftcont} says that the survival
probability is continuous whenever $\lambda_1 \neq \lambda_a$,
thus Conjecture~\ref{cj2} implies Conjecture~\ref{cj3} for trees
of strongly exponential growth.  The reason some kind
of growth condition is needed in Conjecture~\ref{cj2} is that 
it is easy to make both $\lambda_1$ and $\lambda_a$ equal to the
critical value for the integer lattice by making the tree contain
indreasingly long stretches of vertices of degree two.
The theorems in this paper make some progress toward proving these
conjectures, separating $\lambda_1$ and $\lambda_a$, bounding 
$\lambda_a$ and $\lambda_2$ near each other, and proving 
continuity of the survival probability, all in special cases.
Here is a summary.  

Say a tree is homogeneous of degree $n$ if 
every vertex has $n+1$ neighbors.  For homogeneous trees of degree 
$n>2$, Theorem~\ref{eight} states that $0 < \lambda_1 < \lambda_a$.
Additionally, it gives reasonable bounds for the location of
$\lambda_1$, $\lambda_a$ and $\lambda_2$.  So there are at least two 
phase transitions when $n>2$.  

Nonhomogeneous trees come in many 
shapes.  This
paper considers two varieties, namely periodic trees and Galton-Watson
trees.  Since there is no single parameter $n$ measuring the
size of such a tree
(the mean growth of a Galton-Watson tree being an irrelevant statistic
in this setting -- Theorem~\ref{rootnlogn}), bounds on the critical
values are harder to get at.  The main result
here is that $0 = \lambda_1 = \lambda_a = \lambda_b = \lambda_2$
is possible.  The last section deals with the continuity of the 
survival probability as a function of $\lambda$.  This is linked
to the separation of $\lambda_1$ and $\lambda_a$, and to the
question of whether survival is still possible if a single
edge is severed.  The relation is: $\lambda_1 < \lambda_a$ 
$\Rightarrow$ survival with a severed edge whenever survival is 
possible without a severed edge $\Rightarrow$ 
continuity of the survival probability in $\lambda$.  
Thus again the process behaves nicely on homogeneous trees
with $n>2$.

\section{Homogeneous trees} \label{homog}

Fix an integer $n \geq 2$.
Let $\Tree$ be a homogeneous tree where each vertex (also called
a node) has $n+1$ neighbors.
Pick any vertex $\rho$ and call that the root; then each vertex other than
$\rho$ has one neighbor closer to $\rho$ (the parent) and $n$ neighbors 
(children) further from $\rho$.  
Define the contact process on $\Tree$ by letting each
infected node recover at rate 1 and each infected node infect each of its
neighbors at rate $\lambda$ where $\lambda$ is any positive parameter.  Denote 
by $\xi (t)$ the set of infected nodes at time $t \geq 0$.
Thus $\xi$ is viewed as a function from $\R^+$ to $2^\Tree$, or sometimes
as a function from $\R^+ \times \Tree$ to $\{ 0,1 \}$.  
Probabilities referring to initial state $\eta$
at time $0$ will be written as $\P_\eta (\xi (t) \cdots )$.  
When the initial state can be understood, the subscript will be dropped.

\begin{th}
For all $n > 2$, $0 < \lambda_1 < \lambda_2 < \infty$.  For $n = 2$ the
same holds except that the middle inequality may not be strict.
\end{th}
This is a consequence of:
\begin{th} \label{eight}
\begin{eqnarray*}
(i) && {\displaystyle {1 \over n}} \leq \lambda_1 \\[2ex]
(ii) && \lambda_1 \leq {\displaystyle {\sqrt{9 + 16/(n-1)} - 1 
    \over 2(n+1)}} < {\displaystyle {1 \over n-1}} \\[2ex]
&&\mbox{Note that } {\displaystyle {\sqrt{9 + 16/(n-1)} - 1 
    \over 2(n+1)}} \approx  {2 \over 3} 
    \left ( {1 \over n} \right ) + {1 \over 3} \left ( {1 
    \over n-1 } \right ) \\[2ex]
(iii) && {\displaystyle {1 \over 2\sqrt{n}}} \leq \lambda_a \\[2ex]
(iv) && {\displaystyle {4 + 2/(\sqrt{n}-1) -
    \sqrt{ 8 + 16/(\sqrt{n} -1) +4/(\sqrt{n} -1)^2} 
    \over 2(\sqrt{n} -1)}} \leq \lambda_a \\[2ex]
&& \mbox{Note that the LHS is asymptotically ${\displaystyle 
    {4-\sqrt{8} \over 2\sqrt{n}}}$ and is always better 
    than $(iii)$} \\[2ex] 
(v) && \mbox{In particular, } .561722, 
    .425516 , .354248 \leq \lambda_2 \mbox{ when $n$ is respectively }
    2,3,4 \\[2ex] 
(vi) && \lambda_2 \leq {\displaystyle \min ( 2 , {4 \over \sqrt{n} 
    - 4}} ) \mbox{ and } \lim \sup_{n \rightarrow \infty} 
    \lambda_2 \sqrt{n} \leq e = 2.71828 \ldots   \\[2ex]
(vii) && \lambda_c \leq {\displaystyle \max (\lambda_2 , {3 + 
    \sqrt{8n+1} \over 2n - 2}) \mbox{ where asymptotically }
    {3 + \sqrt{8n+1} \over 2n - 2} \approx {\sqrt{2} \over \sqrt{n}}} .
\end{eqnarray*}
\end{th}

The upper bound on $\lambda_1$ in $(ii)$ is less than the lower
bound on $\lambda_a$ on $(iii)$ except for $n \leq 4$.  The reason it
is worth improving the relatively easy bound $(iii)$ to the
messy bound $(iv)$ is to separate $\lambda_1$ and $\lambda_a$
for $n = 3,4$, thus showing there are at least two phase transitions
except possibly when $n = 2$.  The rest of this section is taken up
by the proof of Theorem~\ref{eight}.

\noindent{Proofs:}

$(i)$:  See Liggett (1985) page 166.  Alternatively, see the remark
after the proof of $(ii)$.   $\Cox$

$(ii)$:  Assign  a weight $W(\xi ) = ak(\xi )+bc(\xi )$ to 
each configuration $\xi$ having
finitely many infected nodes, where $a$ and $b$ are constants to be 
determined later, $k$ is the number of infected nodes of $\xi$, and 
$c$ is the number of components in the induced subgraph 
of $\Tree$ on $\xi$, i.e. the graph gotten by restricting $\Tree$ to 
$\xi$.  For $\lambda > (\sqrt{9 + 16/(n-1)} -1)/2(n+1)$, it will be shown 
that $W(\xi )$ is a submartingale, and in fact $\E (W(\xi (t+dt)) \| \xi (t))
\geq (1+\ee ) W(\xi (t))$ for a suitable $\ee > 0$.  Together with the
fact that $W$ is a pure jump process with bounded jumps, this implies
that the probability of $W$ getting arbitrarily large before hitting zero
is bounded below.  (Proof: $1/W$ is a supermartingale when $W$ is greater
than an appropriate function of $\ee$ and the bound on the jump size, so
choose a starting time when $W$ is large enough and stopping when $W$ gets
too small, the optional stopping theorem says the stopping time is 
reached with probability less than $1$.)  Here are some preliminary results.

\noindent{(A)}  The rate of new infections is precisely $\lambda ((n-1)k + 
2c)$.  To see this, count all the attempted infections and subtract off the
times a node attempts to infect an already infected node.  The attempted 
infection rate is $\lambda (n+1)k$ and the failure rate is just $\lambda$
times twice the number of edges having infected nodes on both ends.  But
the number of such edges is the number of infected nodes minus the number
of components they constitute, so the infection rate is $\lambda ((n+1)k
- 2(k - c)) = \lambda (n-1)k + 2c$.

\noindent{(B)}  The total recovery rate is precisely $k$.

\noindent{(C)}  The rate at which new components are formed is precisely 
$k - 2c$, if we agree to count the component lost by the recovery of
an isolated infection as $-1$ new components.  With this convention, 
a new component can only be formed by a recovery, not by a new infection,
and conversely, components can only be lost by the joining 
of previously formed
components along new infections.  Consider any component of $r$ 
infected nodes.  It contains $r-1$ connecting edges.  A recovery at
a node with $j$ edges incident creates $j-1$ new components (this covers
the trivial case $j=0$ as well).  Summing over all nodes gives
$\sum (j-1) = (\sum j) - r = 2r - 2 - r = r-2$, where $\sum j = 2r - 2$ 
because each edge is counted exactly twice.

\noindent{(D)}  The rate at which components are lost by joining is
at most $(n+1)\lambda (c - 1) < 2 \lambda c$.  To see this, consider  
a new, acyclic, bipartite graph, $H$, whose vertices are the components 
of the infection together with all {\em joining vertices}, i.e. vertices 
of the old tree that would join two components by becoming infected.  
The edges of $H$ are all edges 
of the tree that connect a component of the infection to a joining vertex. 
Now, each joining vertex $v$, if infected, causes the loss of 
$e(v)-1$ components, where $e(v)$ is the number of edges of $H$ incident
to $v$.  Since each joining vertex is infected at a rate
of at most $(n+1)\lambda$, summing over all joining vertices gives
a total component loss rate of at most $(n+1)\lambda \, 
\sum_{v \in H} (e(v ) -1)$ which is just $n+1$ times the total number of 
edge of $H$ minus the total number of joining vertices, since 
each edge in the bipartite graph $H$ gets counted exactly once.  
But the total number of edges in the acyclic graph $H$
is just the number of vertices of $H$ minus the number of 
components, which is at most $c$ plus the number of joining vertices 
minus $1$.  The number of joining vertices cancels out and the calculation
is done.

Now fix any $\lambda$, and calculate whether the
infinitesimal expected change in $W$ is always positive.
$$\E (W (\xi (t +dt)) \| \xi (t)) \geq W(\xi (t)) + \left \{ [((n-1)\lambda -1)k
  + 2\lambda c ] a + [k - 2c - (n+1) \lambda ) c] b \right \} dt . $$
This is linear in $c$ and $k$, so to check if it is positive for fixed
$\lambda , a , b$ and every $c$ and $k$, it suffices to check the extreme 
values for $c$, namely $c = 1$ and $c = k$.  In fact it is easier to check
$c = 0$ instead of $c = 1$.  Then the two inequalities to be satisfied
are 
\begin{eqnarray}
((n - 1)\lambda - 1) a + b  & > & 0 \label{c=0} \\[2ex]
((n + 1)\lambda - 1) a - (1 + (n+1)\lambda )b  & > & 0 \label{c=n} .
\end{eqnarray}
By $(i)$, it is safe to assume $(n+1)\lambda - 1 > 0$ and hence to change
equation~(\ref{c=n}) into $a > ((n+1)\lambda +1) b / ((n+1) \lambda - 1)$.
Plugging this into~(\ref{c=0}) shows there will be a solution when
$${\displaystyle {((n - 1)\lambda - 1)((n+1) \lambda + 1) \over 
((n + 1)\lambda - 1)}} + 1 > 0$$
which happens iff $(n+1)\lambda^2 + \lambda - 2/(n-1) > 0$ which, for positive
$\lambda$, will be true when $\lambda > (\sqrt{9 + 16/(n-1)} - 1)/2(n+1)$.  

Finally, it is easy to see that when this inequality holds, 
$[\E (W (\xi (t+dt)) \| \xi (t)) - W(\xi (t)) ] / dt$ is at least 
a small constant multiple of $bn$ for each fixed $\lambda$, thus at 
least $\ee W(\xi (t))$ for an appropriate $\ee > 0$ and the argument 
is done.   $\Cox$

\noindent{Remark:}  This method also gives a proof of $(i)$.  Instead of
the upper bound on component loss in preliminary result~$(D)$, 
use a lower bound of zero.  The same calculation now gives $a = 
b/((n+1)\lambda - 1)$ and plugging that into equation~(\ref{c=0})
with the inequality reversed gives $((n-1)\lambda -1)/((n+1)\lambda -1) 
+ 1 > 0$ which is satisfiable when $\lambda < 1/n$.  Thus for $\lambda < 1/n$,
$W$ can be defined to be a supermartingale and the weight must converge
to zero.

$(iii)$:  This time use a weighting on all the vertices.  It will be shown
that the total weight of the vertices goes to zero even if their cardinality
does not.  Let the weight of a vertex $v$ be $W(v) = n^{-d/2}$ where $d$
is the distance from $\rho$ to $v$, i.e. the
number of edges in the shortest path from $\rho$ to $v$.  
Let 
$$M(t) = \sum_v W(v) I_{v \in \xi (t)} = \mbox{ total weight at time }
    t.$$
Then
\begin{eqnarray*} 
&& \E (M(t + dt) \| \xi (t)) \\[2ex]
& = & \E (\sum_v I(v \in \xi (t+dt)) W(v) \| \xi (t))  \\[2ex]
& \leq & \sum_v I(v \in \xi (t+dt)) W(v) (-1 + \lambda n^{1/2} 
    + n\lambda n^{-1/2}) \\[2ex]
& = & M(t) + dt \, M(t) (2 \sqrt{n}
    \lambda - 1) .
\end{eqnarray*}
So for $\lambda < 1/2\sqrt{n}$, $\E M(t)$ goes to zero exponentially 
fast, hence so does $\P (\rho \in \xi (t))$.   $\Cox$

$(iv)$ and $(v)$:  This time the idea is to take into account the 
coalescing nature
of the infection in order to improve on the bound in $(iii)$.  The
weight of a vertex will depend on whether its parent is infected and will
count less if so, due to the redundancy.  Define
\begin{equation} \label{weight}
W(v) =  r^k ( 1 - d \,I (\mbox{parent of }v \in \xi (t)) ) .
\end{equation}
where $k$ is the distance from $v$ to the root and $r$ and $d$ are
parameters that will be chosen optimally later.  
The goal is to show that the expected change in the total weight $\sum_\Tree
W(v ) I(v \in \xi (t))$ is always negative.  Let $\parent (v )$
and $\child (v )$ denote respectively the number of infected
parents of $v$ (zero or one) and the number of infected children of
$v$.  Define
$$ u(v ) = -1 + \parent (v ) d + \child (v )rd + \lambda [ r^{-1} -
  d + (n - \child (v )) r (1-d)] .$$
Then $u (v ) r^k \, \partial t$ is an upper bound on the
expected increment in the total weight from 
$t$ to $t+\partial t$ due to the recovery of $v$ or infection of other nodes
by $v$.  (Breaks down as follows: recovery loses weight $(1 - \parent (v ))
d r^k$, but gains $\child (v ) d r^{k+1}$ ; infecting the parent gains
at most $r^{k-1}$ while losing $d r^k$, and infecting a child gains
$(1-d)r^{k+1}$ while it may lose something as well.)  

One way to show that the expected increment in total weight is negative
show that $u(v ) < 0$ for each infected node.  Trying
to count this way is inefficient because $d$ must be kept small so that
a node that is completely surrounded by infection does not seem to increase
the total weight if it recovers.  A solution is to
let $Y$ be another parameter to be determined later, and define 
\begin{equation} \label{WEIGHT}
U(v) =  u(v ) + \parent (v ) Y / r - \child (v ) Y .	
\end{equation}
The sum over all nodes of $U (v )r^k$ is the same as the sum of $u(v )r^k$,
because for each parent and child node that are both infected, a 
quantity of $r^k Y$ has merely been transferred from parent to child.
The $U$'s are calibrated so as no longer to penalize our strategy at nodes
with many infected neighbors.

Now there are $2n+2$ cases to check in the forward equation, 
one for each value of $\parent (v )$ and $\child (v )$.
For each pair of these values, $r,d, \lambda$ and $Y$ must be chosen
to make $U < 0$.  Again by linearity, it suffices to consider the cases
where $\child (v ) = 0$ or $n$, so the inequalities are:
\begin{eqnarray*}
-1 + \lambda [ r^{-1} - d + n r (1-d) ] & < & 0  \\[2ex]
-1 + nrd + \lambda ( r^{-1} - d ) - nY & < & 0  \\[2ex]
-1 + d + \lambda n r (1-d) + Yr^{-1} & < & 0  \\[2ex]
-1 + d + nrd + Yr^{-1} - nY & < & 0  .
\end{eqnarray*}

These inequalities are not too hard to solve.  Given $n, r, $ and $\lambda$,
equality in the first equation gives an optimal choice for $d$ and
equality in the third equation gives an optimal choice for $Y$.  It
then turns out that equality in both of the remaining equations can
be achieved for the greatest $\lambda$ when $r = 1/\sqrt{n}$, and hence
the strict inequalities are satisfied for any smaller $\lambda$.  
Plugging in these values
$$r = 1/\sqrt{n} \hspace{.5in} d = {2\sqrt{n} - \lambda^{-1} \over 1 + 
    \sqrt{n}} \hspace{.5in} Y = (1 - \sqrt{n} \lambda) \left ( {1 - 
    \sqrt{n} + \lambda^{-1} \over 1 + \sqrt{n} } \right ) / \sqrt{n}$$
gives a nasty looking quadratic in
$\lambda$ that simplifies to 
$$-(1 - \sqrt{n})^2 \lambda^2 + (4\sqrt{n} - 2)\lambda -2 < 0$$
and has the solution appearing in $(iv)$ (divide through by $1 -
\sqrt{n}$ and use the quadratic formula).  Values for $\lambda , r, d$
and $Y$ for $n = 2,3,4$ are given below so that skeptics may plug them
into the four equations above and check that they actually work.  
For $n = 2$ use $\lambda = .561722$, $r = .7071 , 
d = .434212$ and $Y = .082293$; for $n = 3$ use $\lambda = .425516$, 
$r = .5773 , d = .407781$ and $Y = .089933$; for $n = 4$ use 
$\lambda = .354246$, $r = .5012 , d = .391747$ and $Y = .088349$.

When all the inequalities are satisfied, it is easy to see,
as in the proof of $(ii)$, that in fact $\E (W(\xi (t+dt)) \| \xi (t))
\leq (1 - \ee ) W(\xi (t))$ for a suitable $\ee > 0$.  Then 
again the expected total weight decays
exponentially, hence the probability of finding an infection in 
any finite set decays exponentially and the proof is done.   $\Cox$

$(vi)$:  The $2$ comes from the fact that the critical value for
$\lambda$ on the one-dimensional integer lattice is know to
be at most $2$ and is an upper bound on $\lambda_2$ by 
monotonicity.  For the nontrivial bound, first record two lemmas.
\begin{lem} \label{bin}
Let $M$ be an positive integer valued random variable and pick
$p < \E M$ (with $\E M$ possibly being infinite).  For any $x > 0$,
let $M_x$ be a binomial$(M,x)$ random variable.  Then there is
an $\ee > 0$ for which $\P (M_x \geq 1) > px \wedge \ee$.
\end{lem}

\noindent{Proof:}  Pass to the bounded case as follows: let L be 
large enough so that $\E (M \wedge L) > p$; let $M' = 
M \wedge L$ and $M_x' = \mbox{binomial}(M' , x)$; then $\P (
M_x \geq 1) \geq \P (M_x \geq 1)$ and it suffices to  prove the lemma
for $M'$.  Now
$$\P (M_x' \geq 2) \leq \P (\mbox{binomial}(L,x) \geq 2) =
1 - (1-x)^L - Lx(1-x)^{L-1} = O(x^2) . $$
Also, $\E M_x' \leq \P (M_x' = 1) + L \P (M_x' \geq 2)$ and thus
$$\P (M_x' \geq 1) \geq \P (M_x' = 1) \geq \E M_x' - L \P (M_x' 
\geq 2) > px$$
for small enough $x$ which proves the lemma.   $\Cox$
\begin{lem} \label{recur}
Let $H$ be any nondecreasing function on the nonnegative reals
with $H(x) \geq x$ on some neighborhood of 0.  Suppose $f$ is 
a function on the nonnegative reals that satisfies 
$$\begin{array}{lcrclr}
(1) &&\inf_{0 \leq t \leq L} f(t) & > & 0 &\mbox{   and} \\
(2) && f(t) & \geq & H \left ( \inf_{0 \leq s \leq t-L} f(s) \right ) &
    \mbox{   for } t \geq L 
\end{array} $$ 
for some $L > 0$.
Then $\liminf_{t \rightarrow \infty} f(t) > 0$.
\end{lem}

\noindent{Proof:}  For any $t_0$ and any $\ee > 0$, (2) implies there is 
a sequence $t_1 , t_2 , \ldots , t_k$ with $f(t_i) \geq H(f(t_{i-1})) -\ee 
2^{-i} \geq f(t_{i-1}) - \ee 2^{-i}$ and $t_k \in [0,L]$.  Then $f(t_0)
\geq f(t_k) - \ee$ which is bounded away from $0$ for small $\ee$
by~(1).  $\Cox$

Continuing the proof of $(vi)$,
let $f(t) = \P_\dd (\rho \in \xi (t))$ in an altered contact process 
where only $n$ of the root's $n+1$ neighbors is allowed to
infect the root.  The object will be to get
a recursion $f(t) \geq H(f(t-2r))$.  Here $2r$ is a convenient size for a 
time increment to be chosen later so as to make $H$ monotone increasing
with $H(0) = 0$ and $H'(0) >1$.  Then Lemma~\ref{recur} with $L=2r$
will imply that $\liminf_{t \rightarrow \infty} p(t) > 0$.  

To get this recursion break down $f(t)$ into a product of conditional
probabilities as follows.  Let 
$$H_1 (t) = \P ( v \in \xi (t-r) \mbox{  for some child } v 
\mbox{ of } \rho ) $$
and
$$H_2 (t) = \P (\rho \in \xi (t) \|
v \in \xi (t-r) \mbox{  for some child } v \mbox{ of } \rho ) .$$
Then $f(t) = H_1 (t) H_2 (t)$.  Clearly a lower bound on $H_2$ that
does not depend on $t$ is the probability $p (r) = \P_x (y \in 
\xi (r))$ in a system with only two neighboring nodes $x$ and $y$.
To calculate this, let $Q(s) = \P_x (y \in \xi (s) \| x \in \xi (u)
\mbox{ for all } u < s)$.  Then $Q(0) = 0$ and $Q$ satisfies the
differential equation $Q'(s) = -Q(s) + (1 - \lambda ) Q(s)$ and 
solving the equation gives
\begin{equation} \label{Qdef}
Q(s) = {\lambda \over 1 + \lambda} (1 - e^{-(1+\lambda ) s}) .
\end{equation}
Obtain from this a crude lower bound on $p (r)$ by requiring $x$
to stay infect up to time $r$:
\begin{equation} \label{crude}
H_2 \geq p(r)  \geq  Q(r) \P_x (x \in \xi (s) \mbox{ for all }
    s < r) e^{-r} Q(r) . 
\end{equation}
A better bound is given by conditioning on the recovery time $s$ of 
$x$, requiring $y$ to be infected at time $s$ and remain 
infected through time $r$:
\begin{eqnarray*}
H_2 &\geq & p (r) \\[2ex]
& \geq & Q(r) e^{-r} + \int_0^r \P_x (x \mbox{ is infected to
    precisely time }s ) Q(s) e^{s-r} \, ds \\[2ex]
& = & e^{-r} \left [ Q(r) + \int_0^r Q(s) \, ds \right ] \\[2ex]
& = & {\lambda \over 1 + \lambda } e^{-r} \left [ 1 - e^{-(1 + 
    \lambda ) r} + r - {1 \over 1 + \lambda} (1 - e^{ - (1 +
    \lambda ) r} ) \right ] .
\end{eqnarray*}

Now for a lower bound on $H_1$, let
$M$ be the number of children of $\rho$ that are infected at time $r$,
not counting the one child that is not allowed to infect $\rho$.
The previous calculation says essentially that
each child is infected at time 
$r$ with probability at least $p(r)$, so $\E M \geq n p(r)$.
If each of these nodes is not allowed to be infected by its parent
between time $r$ and time $t-r$ then by definition of $f$ such a 
node infected at time $r$ is infected at time $t-r$ with probability 
$f(t-2r)$.  Thus the actual number of children infected at time
$t-r$ stochastically dominates a random variable that is
$\mbox{binomial}(M,f(t-2r))$.  Then by Lemma~\ref{bin}, $H_1 (t) \geq
n p(r) f(t-2r) \, \wedge \ee$ for some $\ee > 0$ and so 
$f(t) \geq n p(r)^2 f(t-2r) \, \wedge \ee$.  Now writing
$H(x) \geq n p(r)^2 x \, \wedge \ee$ it remains only to show that
$H'(0) > 0$, for then Lemma~\ref{recur} will finish the proof.

For the first claim in $(vi)$ use the crude lower bound~(\ref{crude})
on $p(r)$ to get
\begin{eqnarray*}
\sqrt{H'(0)} & = & \sqrt{n} p(r) \\[2ex]
& \geq  & \sqrt{n} \left ( {\lambda \over 1+ \lambda } \right ) e^{-r}
    (1 - e^{-(1+\lambda )r} ) \\[2ex]
& \geq & \sqrt{n} \left ( {\lambda \over 1+ \lambda } \right ) e^{-r}
    (1 - e^{-r} ) 
\end{eqnarray*}
This holds for any $r$, and an optimal choice of $r$ is $\ln (2)$,
giving $\sqrt{H' (0)} \geq \lambda \sqrt{n} / 4(1+\lambda )$, hence
$H' (0) > 1$ if  $\lambda > 4 / (\sqrt{n} - 4)$.

To get the asymptotic, use the better bound to get
\begin{eqnarray*}
\sqrt{H'(0)} & \geq & \sqrt{n} p(r) \\[2ex]
& \geq & \sqrt{n} \left ( {\lambda \over 1+\lambda } \right ) 
    e^{-r} (1 - e^{-(1+\lambda )r} + r - {1 \over 1+\lambda }
    (1 - e^{-(1+\lambda )r} )) .
\end{eqnarray*}
As $n \rightarrow \infty$ the minimum $\lambda$ satisfying this
inequality will have $\sqrt{n} \lambda \rightarrow e^r \, r^{-1}$
and choosing $r = 1$ gives $\lim \sup \sqrt{n} \lambda = e$.  $\Cox$

$(vii)$:  First an outline of the argument.  By duality, it suffices
to show that when 
$\lambda > \lambda_2$ and $\lambda > (3+\sqrt{8n+1}) / (2n-2)$, 
a process that 
lives forever hits any big set.  By duality again, this means that
the process started from an infection at the root and run for a long
time intersects an independent process initially infected on 
a sufficiently big set and run for a long time.  This will be true
because each of these processes will at some time infect a particular
large configuration if it lives long enough (this uses $\lambda
> \lambda_2$) and the configuration is chosen so that any future
scenarios started from that configuration intersect with high 
probability.  These two facts become the following two lemmas.

\begin{lem} \label{infect.eventually}
Suppose $\lambda > \lambda_2$.  Then for any finite set $S$ and
any initial configuration $\eta$, $\P_\eta (S \subseteq \xi (t)
\mbox{ for arbitrarily large } t) = \P_\eta (\mbox{ infection lives
forever})$. 
\end{lem}

\noindent{Proof:}  Assume without loss of generality that $\eta$
is finite; an easy limit argument handles the infinite case.
Fix $\lambda > \lambda_2$.  The first observation
is that for any node $v$, $\P_v (v \in \xi (t) \mbox{ for arbitrarily 
large } t) = p > 0$,
where $\P_v$ denotes $\P_{\{v\}}$, the contact process run from
only $v$ infected initially.  This is immediate from the definition
of $\lambda_2$, since if $\P_v (\sup \{ t: v \in \xi (t) \} < \infty)
= 1$ then $\P_v (v \in \xi (t) \mbox{ for some } t \geq T) \rightarrow
0$ as $T \rightarrow \infty$ contradicting $\lim\sup_T \P_v
(v \in \xi (T)) > 0$.  Next observe that $\P_\eta (S \subseteq
\xi (t) \mbox{ for arbitrarily large }t \| v \in
\xi (t) \mbox{ for arbitrarily large }t ) = 1$.  This follows from
the Markov property and the fact that $\P (S \subseteq \xi (t+1)
\| v \in \xi (t))$ is bounded away from zero for a fixed
finite $S$, fixed $v \in \Tree$ and any initial state of 
the contact process.  

The last step in proving the lemma is to show that if the process
lives forever then with probability $1$ there is some $v$ that
is infected at an unbounded set of times.  Assume to the contrary.
Then there is some configuration $\eta'$ such that conditional
upon the process being in configuration $\eta'$ at time $t$, the 
probability is greater than $1-p$ that the process lives for ever
but no node is infected at unbounded times.  This contradicts the
definition of $p$ above, since $\eta'$ must have at least one 
infected node, which is then infected at arbitrarily large future
times with probability at least $p$.   $\Cox$
 
\begin{lem} \label{indep.meet}
Define a version of the contact process starting with a single
infection at the root, where no vertex may infect its parent.
Let $\xi (t)$ and $\xi' (t)$ be independent such
processes, started simultaneously.  Then $\P (\xi (t) \cap \xi' (t)
\neq \emptyset \mbox{ for all } t) > 0$, provided $\lambda >
\displaystyle{ {3 + \sqrt{8n+1} \over 2n-2}}$.  
\end{lem}

\noindent{Proof:}  The argument will show that the set of nodes
simultaneously infected in $\xi$ and $\xi'$ contains a supercritical
branching process.  This process $Z$ will be coupled to $\xi$ and
$\xi'$ and for each $v \in \Tree$ the set $\{ t : v \in Z(t) \}$
wil be a possibly empty interval called the lifetine of $v$.
The set of nodes that are ever alive is a branching process, oriented by
the orientation of the rooted tree.  Specifically, let $\rho$
be the single node of $Z$ alive at time zero.  Let each node $v$ of
$Z$ die whenever $v$ recovers in $\xi$ or $\xi'$.
Thus lifetimes in the $Z$ process end at poisson rate $2$.  Dead
nodes may not come back to life.  A node $v$ in the $Z$ process is born
at time $t$ when the following conditions are met.  The parent 
of $v$ must be alive in the $Z$ process, since some time $t_0 < t$. 
The node $v$ must be infected in both $\xi$ and $\xi'$, having just
received an attempted infection from its parent in one of them at 
time $t$ (it is OK if the node was already infected when it 
received the attempt).  Furthermore, if $\xi'$
is the process in which the infection occurs at time $t$, then the
most recent infection of $v$ in the $\xi$ process must have been
after time $t_0$; similarly if the infection at time $t$ occurs in
the $\xi$ process, then the most recent infection of $v$ in the
$\xi'$ process must have been after time $t_0$.  

Whether a node is ever alive in the $Z$ process depends
on: whether the parent was ever alive; the parent's lifetime, 
$t_1 - t_0$; and the attempted infections and recoveries at that 
node in 
the $\xi$ and $\xi'$ processes during the interval $[t_0 , t_1)$.
Since the distribution of attempted infections
in an interval depends on the
interval only through its length, and since lifetimes are i.i.d., 
$Z$ is a Galton-Watson process.  The mean number of offspring
is easy to calculate exactly.  The probability that $v$ is ever alive 
given that the parent is born at time $t$ is just the probability 
that $v$ receives an attempted infection from its parent in both $\xi$
and $\xi'$ before its parent recovers.  This is just $2 \lambda /
(2 \lambda + 2)$ times $\lambda / (\lambda + 2)$, the first of these
being the probability that even one attempt arrives in time, and the
second being the race between the other attempt and the recovery
of the parent.  The branching process is supercritical if $n$ times
this quantity is greater than $1$, which happens when $n \lambda^2 >
\lambda^2 + 3\lambda + 2$, and solving the inequality gives the
condition in the statement of the lemma.   Finally, note that
$Z(t) \subseteq \xi (t) \cap \xi' (t)$, hence $\xi(t) \cap
\xi'(t) \neq \emptyset$ for all $t$ when $Z$ lives forever.   $\Cox$

\begin{lem} \label{fatring}
Let $\eta$ be the configuration where the infected nodes are those
at distance precisely $r$ from the root, for some $r > 0$.  Let
$\lambda$ be as large as in Lemma~\ref{indep.meet}, and let $\xi$
and $\xi'$ be two independent versions of the contact process run
from initial state $\eta$.  Then for any $\ee > 0$, $r$ can 
be chosen large enough so that $\inf_{s,t} \P (\xi (s) \cap \xi' (t) 
\neq \emptyset ) > 1 - \ee$.  
\end{lem}

\noindent{Proof:}  By duality, $\P (\xi (s) \cap \xi' (t) \neq 
\emptyset) = \P (\xi ((s+t)/2) \cap \xi' ((s+t)/2) \neq \emptyset)$,
so it suffices to consider the case $s=t$.
By monotonicity, the processes dominates altered
versions where no node may infect its parent.  These processes proceed
independently on the subtrees below each of the $(n+1) n^{r-1}$ 
nodes in $\eta$, and on each subtree $\xi$ and $\xi'$ are versions
of the process described in Lemma~\ref{indep.meet}.  By independence
and the conclusion of Lemma~\ref{indep.meet}, $r$ can be chosen 
large enough so that $\xi$ and $\xi'$ intersect for all time on
at least one of the subtrees with probability $\geq 1 - \ee$.   $\Cox$

\noindent{Proof of $(vii)$:  } Convergence in distribution is determined
by finite-dimensional marginals, which by inclusion-exclusion
are determined by the probabilities of infection existing in 
each finite set.  So it suffices to show for finite sets $S
\subseteq \Tree$ that $\P_\eta (S \cap \xi (t) \neq \emptyset) 
\rightarrow {\overline \mu} (S \cap \xi \neq \emptyset)$ times 
the survival probability of $\eta$.  By duality, this is just
the product of the survival probabilities of $S$ and $\eta$.  

Let $\ee > 0$ be arbitrary and choose $r$ as in
Lemma~\ref{fatring}.  Let $\xi$ and $\xi'$ be independent versions
of the contact process run from initial states $S$ and $\eta$
respectively.  Let
\begin{eqnarray*}
\tau & = & \inf \{ t \,:\, v \in \xi (t) \mbox{ for all } v
    \mbox{ with } d(\rho , v) = r \} \\
\sigma & = & \inf \{ t \,:\, v \in \xi' (t) \mbox{ for all } v
    \mbox{ with } d(\rho , v) = r \} .
\end{eqnarray*}
Let $\xi''$ and $\xi'''$ be independent contact processes run from
the initial configuration in Lemma~\ref{fatring}, coupled to
$\xi$ and $\xi'$ respectively so that $\xi''(t) \subseteq \xi (t + 
\tau )$ whenever $\tau < \infty$ and 
$\xi''' (t) \subseteq \xi' (t + \sigma)$ whenver $\sigma < \infty$.
Then $\P (S \cap \xi (T) \neq \emptyset ) = \P (\xi' (T/2)
\cap \xi (T/2) \neq \emptyset)$ by duality,
\begin{eqnarray*}
& = & \int \int_{\infty \geq s,t \geq 0}  \P (\sigma \in ds)
   \P (\tau \in dt) \P (\xi' (T/2) \cap \xi (T/2) \neq \emptyset 
   \| \sigma = s , \tau = t) \\[2ex]
& \geq & \int \int_{T/2 \geq s,t \geq 0}  \P (\sigma \in ds)
   \P (\tau \in dt) \P (\xi' (T/2) \cap \xi (T/2) \neq \emptyset 
   \| \sigma = s , \tau = t) \\[2ex]
& \geq & \int \int_{T/2 \geq s,t \geq 0}  \P (\sigma \in ds)
   \P (\tau \in dt) \P (\xi''' (T/2 - \sigma) \cap \xi'' (T/2 - \tau ) 
   \neq \emptyset ) \\[2ex]
& \geq & \int \int_{T/2 \geq s,t \geq 0}  \P (\sigma \in ds)
   \P (\tau \in dt) (1 - \ee ) \\[2ex]
& = & (1 - \ee ) \P (\sigma , \tau \leq T/2) .
\end{eqnarray*}
The last expression converges to $(1 - \ee )$ times the product
of the survival probabilities from $S$ and $\eta$ according to
Lemma~\ref{infect.eventually}, since $\sigma$ and $\tau$ are
independent.  Since $\ee$ is arbitrary, this shows 
$\lim\inf_t \P_\eta (S \cap \xi (t) \neq \emptyset)$ is
at least the product of the survival probabilities.  Duality
gives trivially that $\lim\sup_t \P_\eta (S \cap \xi (t) \neq 
\emptyset)$ is at most the product of the survival probabilities,
and the theorem is proved.   $\Cox$

\section{Nonhomogeneous trees} \label{nonhom}

The upper and lower bounds Theorem~\ref{eight} gives for the critical 
values of $\lambda_2$ on homogeneous trees are within a constant
factor.  Since the constant may be quite large, it was reasonable to
get some more accurate values for small $n$.  In the case of 
nonhomogeneous trees, more work is needed just to get the bounds
within a constant asymptotic factor, as shown by the example in the 
following paragraph.  This section is therefore devoted to improving
the bounds so that they are again within a constant factor of each other.
This section does not discuss bounds on $\lambda_1$ {\em per se}, but
when the trees are nonhomogeneous enough, the better upper bounds on
$\lambda_2$ get pushed down by more than a constant factor
below the upper bounds on $\lambda_1$ that
can be obtained using the methods of the previous section,
which means that the correct value of $\lambda_1$
in these cases is also more than a constant factor lower than can
be gotten by the methods of the previous section.  In fact, the 
last example of this section shows that $\lambda_1$ and $\lambda_2$ may
both be zero on a Galton-Watson tree whose generating function has
all moments.  
For ease of exposition, large constants are chosen with reckless abandon,
leaving the interested reader the exercise of determining better constants.

To see that there is indeed a problem, try to apply the methods of the
previous section to the following tree.  The root has $n$ children, each of 
which has only one child.  Each of the $n$ grandchildren has $n$ children, and
in general the number of children alternates between $1$ and $n$ along
any direct lineage.  

For an upper bound on $\lambda_2$, adapt the method used to 
prove $(iv)$.  Let $f(t) = \P (\rho \in \xi (t))$ as before, and let
$g(t) = \P_{\Tree'} (\rho \in \xi (t))$ where $\Tree'$ is the subtree
rooted at a child of the root of $\Tree$.  Then the same argument
as before gives $f(t) \geq H(g(t-2r))$ and $g(t) \geq J(f(t-2r))$
where $H$ is the same as in the proof of $(iv)$ and $J$ is similar
but based on only $1$ child.  Then monotonicity of $H$ gives
$f(t) \geq (H \circ J)(f(t-4r))$.  Now since $(H \circ J)'(0) = H'(0)J'(0)$,
Lemma~\ref{recur} together with equation~(\ref{crude}) imply that
$\lambda > \lambda_2$ whenever the quantity $Q(r)$ in 
equation~(\ref{Qdef}) satisfies $e^{-4r} Q(r)^4 \cdot
1 \cdot n > 1$.  This means $\lambda$ must be of order $n^{-1/4}$
as $n$ gets large.  Thus if $G$ is the geometric mean of $n$ and $1$,
an upper bound on the second critical value for $\lambda$ is 
\begin{equation} \label{badupper}
\lambda_2 \leq c n^{-1/4} = c \sqrt{G} .
\end{equation}
For trees with longer periodicity, the same argument yields $\lambda_2
\leq c \sqrt{G}$ where $G$ is the geometric mean of the whole sequence
of family sizes.

On the other hand, the somewhat ad hoc methods in $(iii)$ and $(iii')$
can also be adapted to this tree to give lower bounds on $\lambda_2$.
The reader is invited to check
that the best simple weighting scheme along the lines of $(iii)$ for 
trees whose generation sizes alternate between $a$ and $b$ has the weight
of a vertex with $a$ children equal to $\sqrt{a}$ times the weight of each
of its children, and the weight of a vertex with $b$ children equal to
$\sqrt{b}$ times the weight of each of its children (this scheme is not
optimal when the periodicity of the generations is three or more).  Under
this weighting scheme, the total weight goes to zero almost surely 
whenever $\lambda < 1/(\sqrt{a}+ \sqrt{b})$.  When $a=1$ and $b=n$, this
gives the lower bound 
\begin{equation} \label{arithbound}
1/(1+\sqrt{n}) \leq \lambda_a \leq \lambda_2 . 
\end{equation}
Thus the bounds in (\ref{badupper}) and (\ref{arithbound}) have different
asymptotics as $n \rightarrow \infty$, one decreasing as $n^{-1/4}$
and the other as $n^{-1/2}$, in contrast to the case where the tree
was homogeneous and the bounds were always within a constant factor.

Which of these asymptotics for $\lambda_2$ is sharp, if either?  
The somewhat surprising answer is that~(\ref{arithbound}) is sharp,
even though~(\ref{badupper}) is clearly a better
representative for the growth rate of the tree.  In fact for reasonably
regular nonhomogeneous trees, the critical value of $\lambda_2$
is determined by $M = $ the (essential) supremum of the number of children
of a vertex and is at most a constant times $r M^{-1/2}$ where 
$r$ is a logarithmic measure of how far apart vertices with 
$M$ children are from each other.  

There are several kinds of nonhomogeneous trees one might wish to study.
The periodic trees mentioned earlier in this section are one 
such class and Galton-Watson trees are another.  Theorem~\ref{star}  
below will be used to obtain results about both of these cases.  Since the 
proofs for Galton-Watson trees are harder, these will be given first
and a brief version of the other will follow.

Let $\Tree$ be a Galton-Watson tree with generating function $f = 
\sum_{n \geq 0} a_n x^n$, i.e. there
is a root $\rho$ and each vertex $v$ has a random number $C(v )$ 
of children, where $C(v )$ are independent and each is equal to 
$n$ with probability $a_n$.  Assume $f'(1) > 1$ so that the tree is infinite
with probability $p > 0$.  To run the contact process on a Galton-Watson
tree that begins with only the root infected, 
choose $\Tree (\omega )$ according to the Galton-Watson distribution
and let $\xi (t , \omega \| \Tree (\omega ))$ be the
contact process run on the tree $\Tree (\omega )$ starting with only
the root infected.    
\begin{pr} \label{01treelaw}
For a given offspring generating function $f$, let 
$$ q (\Tree ) = \liminf_{t \rightarrow \infty} \P (\rho \in \xi (t) 
\| \Tree (\omega ) = \Tree ) .$$
Then either $q$ is almost surely zero or else 
$\P (q > 0 \| |\Tree | = \infty) = 1$.  
\end{pr}

\noindent{Proof:}  Let $p_1 = \P (q = 0)$.  If the root has $n$ children 
with subtrees $\Tree_1 , \ldots , \Tree_n$ then the $\Tree_i$ are 
$\iid$ with distribution independent of $n$.  Furthermore, if $q(\Tree_i)
> 0$ for any $i$ then $q(\Tree ) > 0$.  To see this use the same
argument as in $(vi)$ of the last section for any fixed $r$ to get 
$$\P (\rho \in \xi (t)) \geq e^{-2r} Q^2 \P_{\Tree_i} (\rho \in
\xi (t - 2r)) . $$
Taking $\liminf$'s over $t$ gives $q(\Tree ) \geq$ a constant
times $q(\Tree_i )$ for
each $i$.  Now by independence of the subtrees $\Tree_i$, 
$p_1  = \P (q(\Tree ) \neq 0) \leq \sum_n a_n (p_1)^n 
= f(p_1 )$.  That
means that either $p_1 = 1$ or $p_1$ is bounded above by the unique
fixed point of $f$ in $(0,1)$.  In the first case $q$ is almost surely 
zero and in the second case $\P (q = 0)$ is at most $\P (|\Tree | < 
\infty )$; since it is also at least $\P (|\Tree | < \infty )$, these 
two must be equal and the proposition is proved.   $\Cox$

In light of this, it is natural to define $\lambda_2 (f)$ to be the $\inf$ of 
$\lambda$ such that $q$ is almost surely nonzero for the contact process with
parameter $\lambda$ run on a Galton-Watson tree with generating function
$f$ conditioned to be infinite.  The following upper bound on 
$\lambda_2$ can now be obtained.
\begin{th} \label{rootnlogn}
Let $r$ be the maximum of $2$ and $c_2 \ln (1/ n a_n) / \ln ( f'(1))$.
There are constants $c_2$ and $c_3$ such that for any $n > 1$, 
$\lambda_2 \leq c_3 \sqrt{r \ln (r) \ln (n) / n}$.
\end{th}

The class of periodic trees may be defined in any of a number of ways
so that a corresponding theorem is true.  Here is one possible course
broad enough to include the example at the beginning of this section.
Let $G$ be any finite graph with a root, $\rho$, a second distinguished
vertex, $\sigma$, and another set of distinguished vertices
$v_1 , \ldots v_m$ distinct from each other.  Suppose that for some 
$j_1$ and $j_2$ a path from $\sigma$ to each $v_i$ of length at most $j_2$
exists that is disjoint from a fixed path of length $j_1$
connecting $\sigma$ and $\rho$.  Let $n$ be $1$
less than the number of neighbors of $\sigma$ and let $j = j_1 + j_2$.  
Construct a graph that will be a $m$-ary tree of $G$'s
as follows. 

Let $\Tree^{(1)} (G) = G$.  For each $i \geq 1$ let $\Tree^{(i+1)} (G)$
be a copy of $G$ together with $m$ copies of $\Tree^{(i)} (G)$ but 
with $\rho$ in the $i^{th}$ copy of $\Tree^{(i)}$ identified with 
$v_i$ in the copy of $G$, so there remains only a single root
in $\Tree^{(i+1)}$, namely the one from the copy of $G$.  The 
increasing limit of these trees (as seen from the root) is the
desired graph $\Tree (G)$.  If $G$ is a tree then so is $\Tree (G)$.
\begin{th} \label{period}
For $G , n , j , m$ and $\Tree (G)$ as above, let $r$ be the maximum of
$2$ and $\lceil j / \ln (m)
\rceil$.  Then there is a constant $c_4$ such that $\lambda_2 \leq
c_4 \sqrt{r \ln (r) \ln (n) / n}$.
\end{th}
The proofs of Theorems~\ref{rootnlogn} and~\ref{period}, 
will be given in section~\ref{pfsegs} since it relies on the results
from section~\ref{finite} on finite stars.

\section{Finite stars} \label{finite}

Let $\Tree$ be a star of size $n$, i.e. $\Tree$ consists of a root, $\rho$
and $n$ other vertices each connected only to $\rho$.  Although any
finite system is eventually trapped in the state of zero infection,
this star turns out to be capable of storing the infection for a 
long time when $\lambda$ is an appropriate constant times 1 / 
$\sqrt{n}$.  Since the application of these results is to 
get the two asymptotics is equations~(\ref{badupper}) 
and~(\ref{arithbound}) to agree within a constant factor, there
is no loss of generality in assuming $n$ is greater than 64, say, 
so that $n^{-1/2}$ and $n^{-1/4}$ are actually distinguishable.
It will also be convenient to restrict $\lambda$ to being less
than $1$; again, nothing is lost, since all the critical values
are less than $1$ for $n \geq 64$.

Before stating the theorems precisely, here is a heuristic explanation
of the results.  Suppose $n$ is large and let $\lambda = a / \sqrt{n}$
for some constant $a$.  Let $x$ be the number of
infected neighbors of $\rho$ divided by $\sqrt{n}$ and pretend 
that $x$ is a continuous rather than a discrete parameter.  
When the root is not infected, $x$ decreases deterministically at
a rate $x$ (i.e. $x = e^{-t} x_0$) and when the root is infected,
$x$ approaches $a$ at rate 1 (i.e. $x - a$ = $e^{-t} (x_0 - a)$).
The root, when infected, recovers at rate $1$, and when healthy,
is infected at rate $(a/\sqrt{n}) (x\sqrt{n}) = ax$.  In other 
words, the state of infection is a Markov process on $\{ 0 , 1 \}
\times [0,a]$ where the second coordinate is $x$ and the first
coordinate is $0$ if the root is healthy and $1$ if the root is
infected.  Its natural scale can be calculated exactly, giving
a function $W : \{ 0 , 1 \}\times [0,a ] \rightarrow \R$
for which $W (t)$ is a martingale.  Writing $f(x) = W(1,x)$ and
$g(x) = W (0,x)$ yields the equations
\begin{eqnarray*}
(a-x) f' + (g-f) & = & 0 \\[2ex]
-xg' + ax (f-g)  & = & 0
\end{eqnarray*}
which have the unique (up to affine transformation) solution
\begin{equation} \label{contmart}
f(x) = \int_0^x e^{-as} (a - s)^{-2} \, ds ~~~; \hfill
g(x) = f(x) - e^{-ax} (a-x)^{-1} .
\end{equation}
From this it is easy to read off the behaviour of $x(t)$ and in
particular, to see that for $k < l < m$, the chance of
hitting $x=a/k$ before hitting $x=a/m$, starting from
$x = a/l$ is on the order $e^{-a^2 (1/l - 1/m)}$.  

Returning to the actual situation, where the number of infections
is integral, there are roundoff errors in the statement of an 
analogous result due to the greatest integer function, and 
roundoff errors in the proof due to the fact that the functions
$f$ and $g$ are no longer very nice.  To kill the roundoff error
I require $n$ to be larger than some $N$ and choose sequentially a 
bunch of constants satisfying certain relations.
The constants involved are named $c_5$ through $c_{13}$ and are
fixed throughout once they are chosen.  These can be evaluated
easily by tracing the argument, but the value of $N$ that this yields
is far from optimal.

\begin{th} \label{star}
Let $c_9 , c_{10} , c_{11}$ be constants satisfying $1/4 >
1/c_{10} + 1/c_{11} > 1/c_{10} > 1/c_9 > 0$ and let $c_5$ 
satisfy $1/c_5 < 1/c_{10} - 1/c_9$.  
For $n$ a positive integer
and $a$ a positive real number, run the contact process on
a star of size $n$ with $\lambda = a / \sqrt{n}$ and initially
$y$ nodes infected, where $y = {\displaystyle \lfloor 
{a\sqrt{n} \over c_{10}} \rfloor}$.  Then the conditions
$$n > (1/c_{10} - 1/c_9 -1/c_5 )^{-1} \hspace{1in} \mbox{ and } 
    \hspace{1in} 4 \leq a \leq \sqrt{n}$$
imply that the probability of the number of infected nodes dropping
below $a\sqrt{n} / c_9$ before reaching $y+a\sqrt{n} / c_{11}$
is at most $e^{a^2 / c_5}$.  
\end{th}

\begin{cor} \label{domingeom}
If $n$ is large enough then
beginning with only $\rho$ infected at time zero, the probability 
of the number of infected nodes staying above $a \sqrt{n} / c_9$ 
during the time
interval $1 \leq t \leq e^{a^2/c_5}/2c_{11}$ is at least $e^{-1}/5$.  
\end{cor}

\noindent{Remarks:}  (1) The important fact here is that the infection 
lives for a time exponential in $a^2$ with a probability bounded away
from zero; this is what creates the order $\sqrt{\ln (n) / n}$ bound
in Theorems~\ref{rootnlogn} and~\ref{period}, via 
equation~(\ref{choosea}) in the next section.  (2) As observed earlier,
requiring $n$ to be large and $\lambda = a/\sqrt{n}$ to be less than
$1$ gives away nothing.

\noindent{Proof of corrolary}:  The probability that $\rho$ remains infected until at
least time $1$ is $e^{-1}$.  Remember this value: $p_1 = e^{-1}$.
Conditioned on that, the probability of
finding at least $y$ nodes infected at time $1$ is at
least the probability of a binomial with parameters $n$ and $Q(1)$ being
least $y$, where $Q(1) = (a / (a+\sqrt{n})) (1-e^{-(1+a/\sqrt{n})})
> (1-e^{-1}) a / 2\sqrt{n}$ is as in equation~(\ref{Qdef}).
The mean of this binomial is at least $a\sqrt{n} /4$ and the 
variance is at most $n Q(1) \leq a n / (a+\sqrt{n})$
so if $n$ is large enough, the binomial is at least $y = a\sqrt{n} /c_{10}$
with probability arbitrarily close to 1, say $4/5$.
Remember this too: $p_2 = 4/5$.
Next, conditional on all that, Theorem~\ref{star} implies 
that the number of times the number of neighbors of $\rho$ goes from 
$y$ to $y+a\sqrt{n}/c_{11}$ before becoming less than $a \sqrt{n} /c_9$ dominates
a geometric random variable with mean at least $e^{a^2/c_5}$.  Finally,
if $n$ is large enough to apply Theorem~\ref{star} then
the waiting times for each new infection are independent from 
everything else and are at least exponential with mean
$1/a\sqrt{n}$ (the infection rate summed over all neighbors of the root
is at most $a\sqrt{n}$).  The mean sum of $a \sqrt{n} / c_{11}$ of these
waiting times $e^{a^2 / c_5}$ times over is thus at least
$e^{a^2 / c_5} / c_{11}$.  Any sum of a geometric random number 
of i.i.d. exponential random variables
is at least half its mean with probability at least $1/4$ 
if the geometric has large enough mean.  
So the probability that the sum of the
waiting times is at least $e^{a^2 / c_5} / 2c_{11}$ conditional
on there being at least $e^{a^2 / c_5}$ transitions from $y$ 
to $y+a\sqrt{n}/c_{11}$ is at least $1/4$.  Call this $p_3$.
Now multiply $p_1$, $p_2$ and $p_3$ to give the result.   $\Cox$

\noindent{Proof} of Theorem~\ref{star}:  For each state $\xi$, 
let $I (\xi )$ be $1$ if the root is infected and $0$ otherwise.
Assign each state $\xi$ the weight
$$W(\xi ) = e^{-ax / 10\sqrt{n}} \left ( 1 - I(\xi )
{\sqrt{n} (e^{a/10\sqrt{n}} -1) \over a}\right ) \;\;, $$ 
where $x$ is the number of 
infected neighbors of $\rho$.  This $W$ is not a martingale like
the heuristic one was, but reducing the exponent from $ax/\sqrt{n}$
to $ax/10\sqrt{n}$ makes it a supermartingale by, as will be evident,
a comfortable margin.  
Note that by the size assumption on $a$, $a/10 \sqrt{n} \leq .1$
so $e^{a/10\sqrt{n}} - 1 \leq 1.2 a/10\sqrt{n}$, hence $\sqrt{n}
(e^{a/10\sqrt{n}} - 1)/a \leq .12$ and $W$ is always positive.
Let $\tau$ be the first time that $\rho$ has 
at least $y+a\sqrt{n}/c_{11} < a\sqrt{n} / 4$ infected neighbors.  
The state $\xi (t)$ is a continuous time Markov chain on the
space $\{ 0 , 1 \} \times \{ 0 , 1 , \ldots , n \}$, where
$\xi$ is represented as the pair ($I(\xi ) , $ number of infections).
The transition rates are: 
\begin{eqnarray*}
(0,x) \longrightarrow (0,x-1) & \mbox{  at rate  } & x \\[2ex]
(0,x) \longrightarrow (1,x) & \mbox{  at rate  } & {ax \over \sqrt{n}} \\[2ex]
(1,x) \longrightarrow (0,x) & \mbox{  at rate  } & 1 \\[2ex]
(1,x) \longrightarrow (1,x-1) & \mbox{  at rate  } & x \\[2ex]
(1,x) \longrightarrow (1,x+1) & \mbox{  at rate  } & {a \over \sqrt{n}} (n - x)
\end{eqnarray*}
For $W(\xi (t \wedge \tau ))$ to be a supermartingale, two inequalities
must be satisfied.  The transitions from $(0,x)$ yield the inequality
$${ax \over \sqrt{n}} {\sqrt{n} \over a} \left ( e^{a / 10\sqrt{n}} - 1 \right )
    \geq x \left ( e^{a / 10\sqrt{n}} - 1 \right )$$
which is satisfied with equality.  Setting $k = (\sqrt{n} / a)(e^{a/10\sqrt{n}}
-1)$, the transitions from $(1,x)$ yield
$$k + x (1-k) (e^{a/10\sqrt{n}} - 1) \leq {a \over \sqrt{n}} (n-x)
    (1-k) e^{-a/10\sqrt{n}} (e^{a/10\sqrt{n}} - 1) $$
and dividing out by $1-k$ and by $e^{a/10\sqrt{n}}-1$ gives
\begin{eqnarray*}
{k \over (1-k)(e^{a/10\sqrt{n}} -1 )} + x & \leq & {a \over \sqrt{n}} (n - x)
    e^{-a/10\sqrt{n}} \\[2ex]
x + {\sqrt{n} \over a - \sqrt{n} (e^{a/10\sqrt{n}} - 1)} & \leq & {a \over
    \sqrt{n}} (n-x) e^{-a/10\sqrt{n}}  \\[2ex]
e^{a/10\sqrt{n}} \sqrt{n} x + n {e^{a/10\sqrt{n}} \over a - \sqrt{n}
    (e^{a/10\sqrt{n}} -1)} & \leq & an - ax \\[2ex]
(a + \sqrt{n} e^{a/10\sqrt{n}}) x  & \leq & n \left ( a - 
    {e^{a/10\sqrt{n}} \over a - 
    \sqrt{n} (e^{a/10\sqrt{n}} - 1 )} \right ) .
\end{eqnarray*}
Since $a \leq \sqrt{n}$, the left hand side is at most
$(1+e^{1/10}) \sqrt{n} x$, and since $a \geq 4$, the right hand 
side is at least $an (1 - e^{1/10}/(4- \sqrt{n}) (e^{2/5\sqrt{n} - 1}))
> an/2$ (the denominator on the right hand side has a local minimum at
either endpoint of $[4,\sqrt{n}]$ with the absolute minimum at $4$). 
Thus the inequalities are satisfied as long as $1.11 \sqrt{n} x
< an/2$ which is always satisfied when $x < a / 4\sqrt{n}$, which
is always satisfied up to the stopping time.  This establishes
that $W(\xi (t \wedge \tau ))$ is a supermartingale.  

To finish proving the theorem, suppose that $\rho$ has $y$ infected 
neighbors at time $t$, where $y$ is as in the statement of the theorem.
There is no loss of generality in the worst-case assumption
that $\rho$ is not infected at time $t$.  Then $W(\xi (t)) =
e^{-ay/10\sqrt{n}}$ and by definition of $y$, this is between 
$e^{-a^2/10c_{10}}$ and $e^{-a^2/10c_{10} + a/10\sqrt{n}}$.  
But $W(\xi ) \geq e^{-a^2/10c_9}$ 
the first time the number of infected nodes falls below 
$a \sqrt{n} /c_9$.  Use the fact that 
$W$ is always positive together with the optional stopping
theorem see that $W$ can reach the value $e^{-a^2 / 10c_9}$ with 
probability at most  $exp (-a^2 / 10c_{10} + a / 10\sqrt{n} + 
a^2 / 10c_9)$.  This gives the desired
conclusion as long as $a^2 / 10c_{10}
- a / 10\sqrt{n} - a^2 / 10c_9 > a^2/c_5$.  
Factoring out $a^2 / 10$ yields the equivalent condition
$1/c_{10} - 1/c_9 - 1/c_5 > 1/a\sqrt{n}$ and since $a$ is at 
least $1$, the condition $n > (1/c_{10} - 1/c_9 - 1/c_5 )^{-1}$
finishes the proof.   $\Cox$

One more lemma will be needed for the next section, so record it now.
\begin{lem} \label{holdsout}
Let $n$, $a$ and $\lambda$ be as in Theorem~\ref{star} and run the contact
process on any graph containing the finite tree $\Tree$ where $\Tree$
consists of a star of size $n$ with the root infected at time zero, to which
has been added a single chain of length $r$ of descendants of some
child $v_1$ of $\rho$, denoted $v_2 , \ldots , v_r$.  Then  
\begin{equation} \label{inho}
\P (v_r \mbox{ is infected before time } e^{a^2/c_5} /2c_{11} ) > 
e^{-1} /5 - (1 - c^r \lambda^r )^{(e^{a^2/c_5}) / 2c_{11}r - 2} 
\end{equation}
where $c$ is a constant independent of everything else.  
\end{lem}

\noindent{Proof:}  By monotonicity, it is sufficient to prove the 
statement for a modified contact process in which $v_2$ cannot infect 
$v_{1}$, so the process restricted to the star of size $n$ is unaffected
by what happens at $v_1$ and below.  Also by monotonicity assume
without loss of generality that the graph is actually equal to $\Tree$.
Let $\B_t$ be the event that the number 
of infected nodes of the star at time $t$ is at least $a\sqrt{n} /c_9$,
$\B = \bigcap \{ \B_t : 1 \leq t \leq e^{a^2/c_5}/1c_{11} \}$,
and let $\F_t$ be the $\sigma$-algebra of events up to time $t$.  Then 
$\P (v_r \mbox{ is infected at some time } s \in  (t , t+r ) \|
\F_t, \B_t)$ is bounded below by the probability that $\rho$ becomes 
infected before time $t+1$ given $\F_t , \B_t$, times the product of 
the probabilities that $v_i$ infects
$v_i+1$ before becoming uninfected and this happens within time $1$.
Working out these probabilities gives a lower bound of $(c \lambda )^r$
for a universal constant $c$.  

Now let $L = \lfloor (1+e^{a^2/c_5})/2c_{11}r \rfloor -1$, and
consider a sequence of times $\{t_i \} =  \{ 1 , 1+r , 1+2r , 
1 + Lr \}$.  For each such time $t_i$, 
$\P (v_r \notin \xi (t) \mbox{ for } t_i < t < t_i + r \| \F_{t_i} ,
\B_{t_i} ) \leq 1 - c^r \lambda^r$, and so
$\P (v_r \notin \xi (t) \mbox{ for } 1 \leq t \leq 1 + Lr \| \B )
\leq (1 - c^r \lambda^r)^L$, since
each $\B_{t_i} \in \F_{t_i}$ and $\B \in \bigcap \B_{t_i}$.
Corollary~\ref{domingeom} gives $\P (\B ) > e^{-1} / 5$ and since
$L$ is at least as big as the exponent in the right hand side 
of~(\ref{inho}) the lemma is proved.   $\Cox$

\section{Proofs of Theorems~\protect\ref{rootnlogn} 
and~\protect\ref{period} and examples} \label{pfsegs}

\noindent{Proof} of Theorem~\ref{rootnlogn}:  Assume without loss of
generality that $a_n > 0$ (or else $r$ is infinite).  
Also, there is no loss of generality in assuming $n$ to be 
large enough to apply Corollary~\ref{domingeom}, since $c_2$
and $c_3$ may be taken large enough so that $c_3 \sqrt{r \ln (r)
\ln (n) / n}$ is at least $2$ for any $n$ that is too small, 
and $2$ is an upper bound on all critical values on trees.
By Proposition~\ref{01treelaw},
is suffices to prove that $\liminf_{t \rightarrow \infty} \P (\rho
\in \xi (t)) > 0$ without conditioning on $\Tree$.  Let
$$g(t) = \P (\rho \in \xi (t)) \geq a_n \P (\rho \in \xi (t) \| 
\rho \mbox{ has at least $n$ children}) . $$  
To calculate $g(t)$, condition on the root having at least
$n$ children until further notice.  Let $c$ be the constant in 
Lemma~\ref{holdsout}.  Treat $r$ as a parameter
to be determined later and choose $a = a(r)$ so that
\begin{equation} \label{choosea}
e^{a^2 / c_5} /2c_{11} r > 2 + 10 (c^r (a/\sqrt{n})^r)^{-1} .
\end{equation}
Using the inequality $(1-x)^{1/x} < e^{-1}$, observe that~(\ref{choosea}) 
forces the right hand side of equation~(\ref{inho}) to be
at least $c_8 = e^{-1} / 4 - e^{-10}$
when $\lambda$ is at least $a / \sqrt{n}$.
To satisfy~(\ref{choosea}), it suffices to choose $a > \sqrt{c_5 
\ln (2c_{11} r \cdot 10 
\cdot c_1^{-r} \sqrt{n}^r) } $ which is at most $ c_3 \sqrt{r 
\ln{r} \ln{n}}$ for an appropriate $c_3$ as long as $r \geq 2$.
Set $\lambda = a / \sqrt{n}$.  If this is sufficient to force 
$\P (q (\Tree ) > 0 \| |\Tree | = \infty ) = 1$ then
$\lambda_2  (f) \leq c_3 \sqrt{ r \ln (r) \ln (n) /n}$
as desired.  Now the root has at least $n$ children; select $n$ of
them and ignore the rest (by monotonicity).
Among all the vertices at distance $r$ from $\rho$,
look for those having at least $n$ children; call these $v_r^1 ,
\ldots , v_r^j$ where $j$ is the random cardinality of the set of
such descendants.  From the facts that $\rho$ has $n$ children, 
that each descendant has an expected $f'(1)$ children and that a node
has $n$ children with probability $a_n$, it follows that
$$\E j \geq n a_n (f'(1))^{r-1} . $$
Now let $M$ be the random number of vertices among $v_r^1 , \ldots , 
v_r^j$ that are infected before $e^{a^2 / c_5} /2c_{11}$.  By the choice of
$a$, Lemma~\ref{holdsout}, and the above expression for $\E j$,
$$\E M \geq c_8 n a_n (f'(1))^{r-1} .$$
Now the argument proceeds as in $(iv)$ of the previous section.  Fix $t$ and
Let $x$ be the minimum value of $g(s)$ for $0 \leq s 
\leq t - e^{a^2 / c_5} /2c_{11}$.  Recursively calculate a lower bound on
$g(t)$ as follows.  

Ignoring all but the first infection of each $v_r^i$ by its parent, 
any of these that are infected at some time $s < t - e^{a^2 / c_5} /2c_{11}$
will evolve independently from time $s$ to time $t - e^{a^2 / c_5} /2c_{11}$
and then be infected with probability at least $x = \inf \{ g(s) \, : \,
s \leq t - e^{a^2 / c_5} /2c_{11} \}$.  Thus the
distribution of the number of such nodes infected at time 
$t - e^{a^2 / c_5} /2c_{11}$ given $M$ dominates a random variable $M_x$ that is
binomial $(M,x)$.  Thus for small enough $x$, $\P (M_x \geq 1) \geq
x (\E M - \ee)$ for arbitrarily small $\ee$.  Finally, if
some $v_r^i$ is infected at time $t - e^{a^2 / c_5} /2c_{11}$, then the
probability of finding $\rho$ infected at time $t$ is bounded below by 
the probability, $p_1$, that the contact process starting with only $v_r^i$
infected at time $t - e^{a^2 / c_5} /2c_{11}$ infects $\rho$ at some time
$s$ with $t - e^{a^2 / c_5} /2c_{11} \leq s \leq t$ times the probability
$p_2$ that $\rho$ is infected at time $t$ given infection of $\rho$
at some such time $s$.  Then $p_1$
is at least $c_8$ by Lemma~\ref{holdsout} and equation~(\ref{choosea})
with the star of size $n$ centered at $v_r^i$ this time.  Meanwhile,
$p_2$ is at least $(e^{-1} / 5) e^{-1} a^2 / (c_9 + a^2)$, which
will be called $c_{13}$.  (See this by applying 
Corollary~\ref{domingeom} together with the fact that the probability
of $\rho$ becoming reinfected between times $t-1$ and $t$ given
at least $a\sqrt{n}/c_9$ infections during this time, and then
staying reinfected, is at least
$a^2 / (c_9 + a^2 )$.)  

Thus as $x \rightarrow 0$, $g(t)$
is bounded below by $c_8 c_{13} x \E (M) \geq 
c_8^2 c_{13} x n a_n (f'(1))^{r-1}$.  Now choose $c_2$ so that putting $r = c_2 \ln 
(1/ n a_n ) / \ln (f'(1))$ forces $c_8^2 c_{13} n a_n (f'(1))^{r-1}$ to be 
greater than $1$.  Applying Lemma~\ref{recur} then gives $\inf_t g(t) > 0$.
$\Cox$

\noindent{Proof} of Theorem~\ref{period}:  Passing to a subtree if
necessary, assume without loss of generality that the parent
of $\sigma$ is on the unique path of length $j_1$ connecting $\rho$
to $\sigma$ and that there is a unique path connecting $\sigma$ to each $v_i$
and it passes through a child of $\sigma$, not through the parent of $\sigma$.

For a positive integer $q$ to be determined later,
consider $\Tree (G)$ as being put together from units
of size $\Tree^{(q)}$ instead of from units of size $G$.
Then $\rho$ has at least $n$ neighbors and furthermore it has at least
$m^{q}$ descendants at distance at most $qj$, all of
which continue recursively the property of having the right number of 
neighbors and well-endowed descendants.  

Now use the argument for the previous theorem with
$m^{q}$ instead of $a_n (f'(1))^{r-1}$.  Then
$a$ must be at least $c_3 \sqrt{qj \ln (qj) \ln (n) / n}$ as before,
and the condition needed for Lemma~\ref{recur} is that $q$ is that $c_8^2 c_{13} 
m^q > 1$.  This means that $q = c' / \ln (m)$ suffices, and setting $ r = 
qj / c'$ proves the theorem.    $\Cox$

A couple of examples conclude this section.  The first shows that the 
the upper bound for $\lambda_2$ in Theorem~\ref{period} can be 
asymptotically sharp.  In fact for this example it is shown that
$\lambda_1$ is within a constant factor of the upper bound for $\lambda_2$
as $n \rightarrow \infty$, so Theorem~\ref{period} provides asymptotically
tight bounds on both critical values.  

\noindent{ {\bf EXAMPLE}} \\
\noindent{This example} calculates $\lambda_1$ for a particular
tree in order to show that the order $\sqrt{\ln (n) / n}$ bound
given by Theorem~\ref{period} may be sharp to within a
constant factor.
Let $\Tree$ be a homogeneous binary tree to every vertex of which 
has been added $n$ children that have no children themselves.  Thus
each vertex has $n+2$ children, only 2 of which have descendants.
To get a lower bound on $\lambda_1$, define the weight of a configuration
as the sum over all nodes that have children of $1$ if the node or
any of its childless children are infected and $0$ otherwise.  The 
following partial converse to Corollary~\ref{domingeom} will be needed.
\begin{pr} \label{converse}
Let $\Tree$ be a star of size $n$ with $\lambda = a / \sqrt{n} < 1$ as
usual.  There is an $N$ for which $n>N$ 
implies that the probability starting from configuration $\eta$ of the
infection dying out before time $\ln (n) + 2 (\ln |\eta | \vee 4)$ 
is at least $e^{-4 a^2} / 10$.
\end{pr}

\noindent{Proof:}  First observe that from any configuration $\eta$,
$\P_\eta (|\xi ( 2 (\ln |\eta | \vee 4))| < 2 a \sqrt{n} ) > 3/4$.  To see
this, couple the process on $\Tree$ to one on another star of size 
$n$, $\Tree'$, where the root of $\Tree'$ always remains infected, but
none of the children of $\rho'$ in $\Tree'$ starts out infected.  The coupling
has corresponding children of $\rho$ and $\rho'$ recover at the same
times and has an infection occurring in $\Tree'$ whenever one occurs 
in $\Tree$ but not necessarily vice versa.  When every child of $\rho$
has recovered once, the number of infected children of $\rho$ is
at most the number of infected children of $\rho'$.  This happens by
time $2 (\ln |\eta | \vee 4)$ with probability $(1 - e^{-2 (\ln |\eta |
\vee 4)})^{|\eta |} > 7/8$.  But the number of infected children of
$\rho'$ at any time is a random variable converging upwards to
a stationary distribution that has mean approximately $a \sqrt{n}$
and for large enough $n$ has probability less than $1/8$ of being 
above $2 a \sqrt{n}$.  Thus the coupling establishes the observation.
Now from a configuration where at most $2 a \sqrt{n}$ nodes are infected,
it remains to show that the probability of the infection dying out 
within time $\ln (n)$ is at least $e^{-4a^2} / 7.5$.
The probability that no reinfection ever occurs is at least $(1 / (
1 + \lambda ))^{4 a \sqrt{n}} > (1 - a / \sqrt{n})^{4 a \sqrt{n}}
\approx e^{-4a^2}$.  If no reinfection occurs, the infection dies
out within time $\ln (n)$ with probability at least $(1 - 
e^{-\ln (n)} )^{2a\sqrt{n}} = (1 - 1/n )^{2a\sqrt{n}} \geq
e^{-2}$ since $a \leq \sqrt{n}$ by assumption. 
Multiplying together the conditional probabilities
gives $e{-4a^2} (3/4) / e^2 \geq e^{-4a^2} / 10$.  $\Cox$

What is needed below is actually an altered version of this proposition
saying that even if $\rho$ is given 3 additional neighbors that are 
always infected, the probability of the infection dying out except on
the additional 3 neighbors in the allotted time is still at
least $e^{-4a^2} / 10$.  The reader may fill in the details to beef up the proof.

Going back now to the example, consider any node $v$ that contributes $1$ to
the weight of the configuration either by having an infected child or
by being a node with children that is itself infected.  There are
only three neighbors of $v$ that may increase the weight of the 
configuration if they get infected by $v$, so the expected contribution
of $v$ to increasing the weight of the configuration in a time period
of length $3 \ln (n)$ is at most 
$$9 \lambda \ln (n) = 9 a \ln (n) / \sqrt{n} . $$
The probability of $v$ decreasing the weight during this time 
because the infection dies out on $v$ and its childless children 
is at least $e^{-4a^2} / 10$, by the beefed up version of the previous
proposition (since $|\eta |$ is necessarily at most $n$).  Now pick
$a = \sqrt{\ln (n)} / 10$.  This makes
$$e^{-4a^2} / 10 = n^{.1} / 10\sqrt{n} / 2 > 9 \sqrt{\ln (n)} 
\ln (n) / 10 \sqrt{n} = 9a \ln (n) / \sqrt{n}$$
for large enough $n$.  Thus any node contributes an expected weight
that is less than $1$ before it and its childless children all recover.
The usual supermartingale argument then implies that the infection dies out
almost surely and so 
\begin{equation} \label{lowerlil}
\lambda_1 \geq \sqrt{\ln (n)} / 10 \sqrt{n}
\end{equation}
for large enough $n$.  

On the other hand, $\Tree$ is a periodic tree: $\Tree = \Tree (G)$ for
the graph $G$ with root $\rho = \sigma$ and $n+2$ other vertices
all connected only to the root, two of which are labelled $v_1$ 
and $v_2$.  In the notation of section~\ref{nonhom}, $m = 2$, $j = 1$
and what is called $n$ in Theorem~\ref{period} is $n+2$.  Then 
Theorem~\ref{period} gives
\begin{equation} \label{upperlil}
\lambda_2 \leq c_4 \sqrt{ \ln {n}} / \sqrt{n} .
\end{equation}
Thus for this class of trees, both~(\ref{lowerlil}) and~(\ref{upperlil})
are sharp to within a constant factor.

\noindent{{\bf EXAMPLE}} \\
\noindent{This} example is thrown in to show that even medium-sized tails on the 
generating function of a Galton-Watson tree can force the critical
values down to zero.  For any positive $\gamma < 1$ let
$$f(x) = \sum_n c e^{-n^\gamma} x^n$$
be the generating function for the number of children, where $c$ 
is a normalizing constant.  Then $\ln (1 / n a_n) \approx 
n^\gamma$, so Theorem~\ref{rootnlogn} gives $\lambda_2 \leq
c\sqrt{n^{\gamma - 1} \ln (n)}$ which goes to zero as $n \rightarrow
\infty$.  Thus $\lambda_2$ is less than any positive number, hence
$\lambda_1 = \lambda_2 = 0$.

\section{Continuity of the survival probability} \label{cont}

Continuity from the right of the survival probability
as a function of $\lambda$ is immediate
from the monotonicity, since versions of the process for each
$\lambda$ can all be coupled so that the event of survival
is continuous from the right in $\lambda$.  At the moment, I can
only prove continuity from the left when the phase transitions
at $\lambda_1$ and at $\lambda_a$ are distinct.  It is also
necessary to assume that the tree is homogeneous or periodic.
In particular, the survival probability is a continuous function
of $\lambda$ for homogeneous trees with $n > 2$.
I believe that survival probability is continuous for arbitrary
trees.  

I will present the argument that $\lambda_1 < \lambda_a$ implies
continuity of the survival probability only for homogeneous
trees; the generalization to periodic trees is straight-forward but
tedious.  In the foregoing discussions, let $p(\lambda )$
denote the survival probability for the contact process on an
$n$-ary homogeneous tree with only the root infected initially.  Begin
with a counting lemma.

\begin{defn}
For a finite set $S \subseteq \Tree$, say a vertex $v \in S$ is 
surrounded in $S$ iff each of the $n+1$ components of $\Tree
\setminus \{ v \}$ intersects $S$.
\end{defn}

\begin{lem} \label{surround}
Let $S$ be a finite set of vertices of a homogenous tree, $\Tree$.  
Then the number
of surrounded vertices of $S$ is less than $|S| / n$.
\end{lem}

\noindent{Proof:}  The proof proceeds by repeated reductio ad
absurdum.  Suppose for a contradiction that the result
is false.  Let $S$ be a counterexample with as few vertices as 
possible, and among such counterexamples, let $S$ have as small
a value as possible for the sum $\sum_{v,w \in S} d(v,w)$.  The
{\em induced subgraph} on $S$ is just the vertices of $S$ together
with whichever edges of $\Tree$ connect two vertices in $S$.  

I claim that the induced subgraph on $S$ must be connected.
Suppose not.  Let $v \in \Tree \setminus S$ be a vertex for which
$\Tree \setminus v$ has at least two components intersecting $S$
(the existence of such a $v$ is equivalent to $S$ failing to be
connected).  There are two cases.  Firstly, suppose more than
two components of $\Tree \setminus v$ intersect $S$.  Writing
$S= \bigcup S_i$, where each $S_i$ is the intersection of $S$
with a single component of $\Tree \setminus v$, deletion of
any single one of the $S_i$ from $S$ does not affect which of the
remaining vertices are surrounded.  Then by deleting the $S_i$ 
with the smallest fraction of surrounded nodes, a smaller 
counterexample is obtained, which is a contradiction.  The other
case is that exactly two components of $\Tree \setminus v$ intersect
$S$.  Let $S_1$ and $S_2$ be the intersection of $S$ with these
two components and let $v_1$ and $v_2$ be the neighbors of $v$
in the two components.  Let $\phi$ be any graph-automorphism of 
$\Tree$ sending $v_2$ to $v$ and $v$ to $v_1$.  Let $S'$ = 
$S_1 \cup \phi[S_2]$.  Then vertices in $S_1$ are surrounded in
$S$ iff they are surrounded in $S'$, and vertices in $S_2$ are
surrounded in $S$ iff their images under $\phi$ are surrounded in $S'$.
So $S'$ is a counterexample with as many vertices as $S$ but a
smaller value $\sum_{v,w \in S} d(v,w)$, which is a contradiction.

Now that it is established that the induced subgraph on
$S$ is connected, count the edges of $S$ in two different 
ways.  First, since the induced subgraph is a tree, there is one
fewer edge than there are vertices.  Counting oriented edges 
doubles this, so there are $2|S| - 2$ oriented edges.  On the other
hand, there are $n+1$ oriented edges leading out of every surrounded
vertex and at least one edge leading out from every vertex, so 
the number of oriented edges is at least $|S| + n $(number of 
surrounded vertices), thus $2|S| - 2 \geq |S| + n $(number of
surrounded vertices), and there are less than $|S| / n$ surrounded
vertices.  This is not a counterexample, thus
no counterexample exists.   $\Cox$

The rest of the discussion deals with the contact process on a 
tree with a severed edge.  For each positive integer $r$, let
$\xi^{(r)}$ be a contact process run from an initial infection
at the root with an edge severed at distance $r$ from the root,
i.e. no infection passes over that edge.  Write $\xi^{(r;\lambda)}$
or $\xi^{(;\lambda)}$ when it is necessary to emphasize the choice of $\lambda$.
An adjacent edge counts as being at distance one.  The main lemma on
continuity from the left is the following.

\begin{lem} \label{severed.cont}
Let $\Tree$ be a homogeneous tree.
For a fixed $\lambda$, let $\xi^{(r)}$ be as above, and suppose the
survival probability for the infection in $\xi^{(r)}$ is positive
for some $r$.  Then the survival probability for the unaltered contact
process is continuous from the left for that value of $\lambda$.
\end{lem}

\noindent{Remark:}  There is a corresponding lemma in percolation
theory stating that the percolation probability is continuous provided
that percolation in $\Z^d$ implies percolation in a half-space.

\noindent{Proof:}  First note that $r$ can be taken to be $1$, since
the initial infection in the $\xi^{(1)}$ process has a nonzero chance
of moving to a distance $r$ from the severed edge before it dies,
thus dominating a $\xi^{(r)}$ process a postive fraction of the time.

The next claim is that if $\xi^{(1)}$ survives with positive probability
at a particular value of $\lambda$ then it survives with positive
probability for some smaller value $\lambda'$.  To find such a 
$\lambda'$ begin by setting $p = p (\lambda)$ equal to 
the survival probability
for $\xi^{(1;\lambda)}$ and pick $t$ large enough
so that 
$$\P (|\xi^{(1;\lambda)} (t)| \geq 10 (1+\lambda ) / \lambda p) \geq 3p/4 . $$
such a $t$ exists, since $|\xi^{(1)} (t)| \rightarrow \infty$ 
almost surely when the infection lives forever.  Pick $\lambda'
\in [\lambda / 2 , \lambda)$ and large enough so that 
$$ \P (|\xi^{(1;\lambda')} (t)| \geq 5 (1+\lambda' ) / \lambda' p) \geq p/2 .$$
This is possible because $p/2 < 3p/4$, $5(1+\lambda' ) / \lambda' p
< 10 (1+ \lambda ) / \lambda p$, and the distribution of 
$\xi^{(1)} (t)$ is continuous in $\lambda$ for fixed $t$.  

For this paragraph, let the infection rate be understood to be
$\lambda'$.  Let $S \subseteq \xi^{(1;\lambda')} (t)$ be the random subset
of the infected nodes at time $t$ that are not surrounded in the
set of infected nodes at time $t$.  Since $1 - 1/n \geq 1/2$, 
Lemma~\ref{surround} shows that $S$ contains at least half the nodes
infected at time $t$.  Thus
\begin{equation} \label{Sisbig}
\P (|S| \geq 5 (1+\lambda' )/ 2\lambda' p ) \geq p/2 . 
\end{equation}
For each $v \in S$ there is by definition of $S$ some neighbor,
say $g (v)$ for which the component of $\Tree \setminus v$ containing
$g(v)$ does not intersect $\xi^{(1;\lambda')} (t)$.  For each such $v$, the
probability is $\lambda' / (1+\lambda' )$ that $v$ will pass the 
infection on to $g(v)$ before recovering.  Let $A$ be the random set 
of such $v$.  Forbid any infection of nodes in $S$ from time $t$
onward.  Then the process proceeds independently on each component
of $\Tree \setminus S$ containing some $x \in A$, since
communication has been severed.  Furthermore, the starting 
configuration on each such component is precisely a single
infection and the process evolves as $\xi^{(1;\lambda')}$, since an edge 
next to the initial edge has been severed.  Then the $\xi^{(1)}$ 
process dominates a set of independent $\xi^{(1)}$ processes, 
started from each vertex in $A$.  As in the theory on Galton-Watson
processes, the survival probability is greater than the fixed
point of a generating function for the probabilities of the various
cardinalities of $A$.  This must be positive if $\E |A| > 1$.
Since $E (|A| \| |S| > 5 (1+\lambda' )/2\lambda p) > 5/2p$, 
the desired result follows from equation~(\ref{Sisbig}).

Finally, to show continuity from the left, recall $p (\lambda' )$
is the survival probability for a severed edge and parameter
$\lambda'$, and let $q$ be the survival probability for the
unaltered contact process with parameter $\lambda$.
Let $\ee > 0$ be arbitrary and let $M$ be large enough so that
$(1-p (\lambda'))^M < \ee / 4$.  Also pick $M$ large enough so that any
binomial random variable with mean at least $M$ is at least $M/2$
with probability at least $1 - \ee / 4$.  Let $t$ be large 
enough so that 
$$\P (|\xi^{(;\lambda)} (t)| \geq 4 (1+\lambda' )M / \lambda') \geq 
   q - \ee /4 $$
for the parameter value $\lambda$, and pick $\lambda''$
with $\lambda' \leq \lambda'' < \lambda$ large enough so that
$$\P (|\xi^{(\lambda'')} (t)| \geq 2 (1+\lambda'' )M / \lambda'' ) 
   \geq q - \ee /2. $$
Rerun the previous argument but for infection rate $\lambda''$.  
Let $S$ be the unsurrounded infected
nodes in $\xi^{(;\lambda'')} (t)$, so each $v \in S$ has a neighbor $g(v)$ 
in a component of $\Tree \setminus S$ with no infection.  There
is a set $A$ of nodes $g(v)$ for $v \in S$, that are infected before 
$v$ recovers.  Formally $A$ is the set of nodes $g(v)$ sucht that
$v \in S$ and there exists a $t' > t$ for which $v$ is infected 
throughout the time interval $[t,t']$ and such that $v$ infects
$g(v)$ at time $t'$.  The binomial hypothesis on $M$ guarantees that
$A$ has cardinality at least $M$ with probability at least 
$q - 3\ee /4$, since each $g(v)$ is independently in $A$ with 
probability $\lambda'' / (1+\lambda'')$.
Each node in $A$ defines an independent process
which is a restarted version of $\xi^{(1)}$ with parameter
$\lambda'' \geq \lambda'$, so the hypothesis on $p(\lambda' )$
and $M$ gives a survival probability of at least $q-\ee$ for
the whole process.  Since $\ee$ was arbitrary, the lemma is proved.  
$\Cox$

It remains to show that the hypothesis of this lemma is satisfied
when the phase transitions at $\lambda_1$ and $\lambda_a$ are
distinct.  This is done in two steps.

\begin{lem} \label{cantreachfar}
Suppose the contact process is run on a homogenous
tree for some $\lambda < \lambda_a$.
Then $\P (v \mbox{ is ever} \mbox{ infected}) \rightarrow 0$ as 
$d(\rho , v) \rightarrow \infty$.
\end{lem}

\noindent{Proof:}  Suppose not.  Then by monotonicity it decreases
to some bound $b> 0$.  Since $\lambda < \lambda_a$ the root
is infected finitely often almost surely, so there is some 
time $t$ and configuration $\eta$ such that $\P (\mbox{survival
without infecting the root any more} \| \xi (t) = \eta) > 1-b$.
Since any node alive at time $t$ infects the root with probability
at least $b$, this is a contradiction.   $\Cox$

\begin{lem} \label{somerworks}
Suppose $\lambda < \lambda_a$.  Then the survival probability
for $\xi^{(r)}$ on a homogenaous tree approaches the survival 
probability for the 
unaltered contact process as $r \rightarrow \infty$.
\end{lem}

\noindent{Proof:}  Couple all the $\xi^{(r)}$ with
the unaltered process, so that they evolve identically until
an infection in the unaltered process crosses an edge that is 
severed in the other process.  By Lemma~\ref{cantreachfar}
the probability of the processes ever becoming uncoupled
goes to zero as $r$ goes to infinity, and the result follows.
$\Cox$

Putting all these lemmas together gives
\begin{th} \label{leftcont}
If $\lambda_1 < \lambda_a$ then the survival probability for
the contact process from a single initial infection on a 
homogeneous tree is continuous in $\lambda$.
\end{th}

\noindent{Proof:}  When the survival probability is zero there
is nothing to prove.  When it is positive and $\lambda < \lambda_a$
then some $\xi^{(r)}$ also has positive survival probability
by Lemma~\ref{somerworks}, and continuity follows from
Lemma~\ref{severed.cont}.  For $\lambda \geq \lambda_a$, the
existence of a $\lambda$ between $\lambda_1$ and $\lambda_a$
implies the survival with positive probability of a severed-edge
contact process for that and hence all higher $\lambda$, and
Lemma~\ref{severed.cont} again applies.  $\Cox$

\renewcommand{\baselinestretch}{1.0}\large

\end{document}